\documentclass{amsart}
\usepackage[utf8x]{inputenc}
\usepackage{amsmath}
\usepackage{amssymb}
\usepackage{psfrag}
\usepackage[all,cmtip]{xy}
\usepackage{amsmath,amscd}
\usepackage{amssymb}
\usepackage{graphicx}
\usepackage{wasysym}
\usepackage{textcomp}

\usepackage{graphicx} %If you want to include postscript graphics
\DeclareGraphicsExtensions{.pdf,.eps,.fig}

\newtheorem{definition}{Definition}
\newtheorem{proposition}{Proposition}
\newtheorem{lemma}{Lemma}
\newtheorem{theorem}{Theorem}

\newtheorem{corollary}{Corollary}
\newtheorem{conjecture}{Conjecture}
\title{Groupoids and Poisson sigma models with boundary
}

\author{Ivan Contreras and Alberto S. Cattaneo}
\thanks{Partially supported by SNF Grant 20-131813.}
\date{May 30 2012}
\begin{document}
\maketitle

\begin{abstract}
This note gives an overview on the construction of symplectic groupoids as reduced phase spaces of Poisson sigma models and its 
generalization in the infinite dimensional setting (before reduction). 
\end{abstract}
%\tableofcontents
\begin{section}{Introduction}
In \cite{Ca}, it was proven that the reduced phase space of the Poisson sigma model under certain boundary conditions and assuming it is a smooth manifold, has the structure 
of a symplectic groupoid and it integrates the cotangent bundle of a given Poisson manifold $M$. This is a particular instance 
of the problem of integration of Lie algebroids, a generalized version of the \textbf{Lie third theorem}\cite{Duistermaat} . 
The general question can be stated as
\begin{itemize}
 \item Is there a Lie groupoid $(G,M)$ such that its infinitesimal version corresponds to a given Lie algebroid $(A,M)?$
\end{itemize}
For the case where $A=T^*M$ and $M$ is a Poisson manifold the answer is not positive in general, as there are topological 
obstructions encoded in what they are called the monodromy groups \cite{Cra}. A Poisson manifold is called integrable if such Lie groupoid $G$ exists.
The properties of $G$ are of special interest in Poisson geometry, since it is possible to equip $G$ with a symplectic structure $\omega$ compatible 
with the multiplication map in such a way that $G$ is a symplectic realization for $(M,\Pi)$.\\
For the integrable case, the symplectic groupoid integrating a given Poisson manifold $(M, \Pi)$ is constructed explicitely in \cite{Ca}, as the phase space modulo 
gauge equivalence of the Poisson Sigma model (PSM), a 2-dimensional field theory.\\

In a more recent perspective (see \cite{AlbertoPavel1, AlbertoPavel2}), the study of the phase space before reduction plays a crucial role. This allows dealing with with nonintegrable Poisson structure,for which the reduced phase space is singular, on an equal footing as the integrable ones. This new approach differs from the stacky perspective of Zhu and Tseng (see \cite{Chenchang}) and seems to be better adapted to symplectic geometry and to quantization. \\

In order to include this construction in the more general setting, where the reduced phase space can be singular and more general source spaces are allowed, the study of the phase space before reduction plays a crucial role.

In a more recent perspective, in order to include this construction in the more general setting, where the reduced phase space can be singular 
and more general source spaces are allowed, the study of the phase space before reduction plays a crucial role.\\
In a paper in preparation \cite{Relational}, we introduce a more general version of a symplectic groupoid, called 
\textit{relational symplectic groupoid}. In the case at hand, it corresponds to an infinite dimensional 
weakly symplectic manifold equipped with structure morphisms (canonical relations, i.e. immersed Lagrangian submanifolds) compatible with the 
Poisson structure of $M$. In this work, we prove that
\begin{enumerate}
 \item  For any Poisson manifold $M$ (integrable or not), the relational symplectic groupoid always exists.
 \item  In the integrable case, the associated relational symplectic groupoid is equipped with locally embedded Lagrangian submanifold.
 \item  Conjecturally, given a regular relational symplectic groupoid $\mathcal{G}$ over $M$ (a particular type of object that admits symplectic reduction), there exists a unique 
Poisson structure $\Pi$ on $M$ such that the symplectic structure $\omega$ on $\mathcal G$ and $\Pi$ are compatible. This is still work in progress.
\end{enumerate}
This paper is an overview of this construction and is organized as follows. Section 2 is a brief introduction to the Poisson sigma model and its reduced phase space.
Section 3 deals with the version before reduction of the phase space and the introduction of the relational symplectic groupoid.
An  interesting issue concerning this construction is the treatment of non integrable Poisson manifolds: even if the reduction does not exists as a smooth manifold, the relational symplectic groupoid always exists. One natural question arising at this point is:
\begin{itemize}
 \item Can there be a finite dimensional relational symplectic groupoid equivalent to the infinite dimensional one for an arbitrary Poisson manifold?
\end{itemize}

The answer to this question is work in progress and it will be treated in a subsequent paper.\\

Another aspect, which will be explored later, is the connection between the relational construction and the Poisson Sigma model with branes, where the boundary conditions are understood as choices of coisotropic submanifolds of the Poisson manifold. The relational symplectic groupoid seems to admit the existence of branes and would explain in full generality the idea of dual pairs in the Poisson sigma model with boundary \cite{Coiso, Branes}.\\
This new program might be useful for quantization as well. Using ideas from geometric quantization, what is expected as the quantization of the relational symplectic groupoid is an algebra with a special element, which fails to be a unit, but whose action is a projector in such a way that on the image of the projector we obtain a true unital algebra. Deformation quantization of a Poisson manifold could be interpreted in this way.

\section{PSM and its reduced phase space.} We consider the following data

%In particular,  the  case when Lie algebroid is $(T^*M,M)$, with $M$ a Poisson manifold has particular interests for Lie theory, since the construction of a smooth groupoid with a compatible symplectic structure gives, between other things, a symplectic realization of the given Poisson manifold. More precisely, it solves the following question:
%\\ \textit{Given a Poisson manifold $(M, \Pi)$, is there a symplectic groupoid $(G,M)$ such that its space of objects is precisely $M$, the (target) source is an (anti) Poisson map and $M$ is a Lagrangian submanifold of G?}\\
%In general a Poisson manifold can fail to be integrable and this means that a groupoid satisfying these properties might be non smooth. The obstruction to integrability of Poisson manifolds is of topological nature \cite{Cra} and is based on the computation of what is called \textit{monodromy groups} of the infintesimal foliation whose space of leaves determines the groupoid $G$.
%In the construction given by Cattaneo and Felder in \cite{Ca}, the symplectic groupoid is constructed as the space of reduced boundary fields of a 2 dimensional topological field theory, called Poisson Sigma model and denoted by PSM, where the source space corresponds to a rectangle with particular boundary conditions and the target space is a given Poisson manifold $M$. More precisely, the data to describe the PSM are:
\begin{enumerate}
\item A compact surface $\Sigma$, possibly with boundary, called the source space.

\item A finite dimensional Poisson manifold $(M,\Pi)$, called the target space. Recall that a bivector field $\Pi \in \Gamma(TM \wedge TM) $ 
is called Poisson if the the bracket $\{,\}: \mathcal{C}^{\infty}(M) \otimes \mathcal{C}^{\infty}(M)\to \mathcal{C}^{\infty}(M)$, defined by
\[\{f,g\}= \Pi(df,dg)\]
is a Lie bracket and it satisfies the Leibniz identity
\[\{f,gh\}=g \{f,h\}+h\{f,g\}, \forall f,g,h \in \mathcal{C}^{\infty}(M).\] 
In local coordinates, the condition of a bivector $\Pi$ to be Poisson reads as follows
\begin{equation}\label{SN}
\Pi^{sr}(x)(\partial_r) \Pi^{lk}(x)+\Pi^{kr}(x)(\partial_r) \Pi^{sl}(x)+\Pi^{lr}(x)(\partial_r) \Pi^{ks}(x)=0,
\end{equation}

that is, the vanishing condition for the Schouten-Nijenhuis bracket of $\Pi.$

\end{enumerate}
The space of fields for this theory is denoted with $\Phi$ and corresponds to the space of vector bundle morphisms between $T\Sigma$ and $T^*M$.
This space can be parametrized by the pair $(X, \eta)$, where $X$ is a $\mathcal C^{k+1}$-map from $\Sigma$ to $M$ and $\eta \in 
\Gamma^k(\Sigma, T^*\Sigma \otimes X^*T^*M).$\\
On $\Phi$, the following first order action is defined:
\[S(X,\eta):= \int_{\Sigma} \eta\wedge dX + \frac 1 2 \Pi^{\#}(X) \eta \wedge \eta,\]
where $\Pi^{\#}$ is the map from $T^*M\to TM$ induced from the Poisson bivector $\Pi$, the integrand, called the Lagrangian, will be denoted by $\mathcal{L}$. Associated to this action, the corresponding variational problem $\delta S=0$
induces the following space 
\[\mbox{EL}=\{\mbox{Solutions of the Euler-Lagrange equations}\}\subset \Phi,\]
where, using integration by parts \[\delta S= \int_{\Sigma} \frac{\delta \mathcal{L}}{\delta X} \delta X+ \frac{\delta \mathcal{L}}{\delta \eta} \delta \eta + \mbox{boundary terms}.\]The partial variations correspond to:
\begin{eqnarray}
\frac{\delta \mathcal{L}}{\delta X}&=& dX+ \Pi^{\#}(X)\eta=0\\
\frac{\delta \mathcal{L}}{\delta \eta}&=& d\eta+ \frac 1 2 \partial \Pi^{\#}(X)\eta\wedge\eta=0.
\end{eqnarray}

Now, if we restrict to the boundary, the general space of boundary fields corresponds to

\[\Phi_{\partial}:=\{\mbox{vector bundle morphisms between  } T (\partial \Sigma) \mbox{  and  } T^*M\}.\]
 
Following \cite{Corfu}, $\Phi_{\partial}$ is endowed with a symplectic form and a surjective submersion $p: \Phi \to \Phi_{\partial}$. We define 
$$L_{\Sigma}:= p(EL)$$.

Finally, we define $C_{\Pi}$ as the set of fields in $\Phi_{\partial}$ which  can be completed to a field in  $L_{\Sigma^{'}}$, with $\Sigma^{'}:= \partial \Sigma \times [0, \varepsilon]$, for some $\varepsilon$.  \\
%the constraint to the boundary is defined as
 It turns out that $\Phi_{\partial}$ can be identified with $T^*(PM)$, the cotangent bundle of the path space on $M$ and that
 \[C_{\Pi}:= \{(X,\eta)\vert dX= \pi^{\#}(X)\eta,\, X: \partial \Sigma \to M, \, \eta \in \Gamma (T^*I \otimes X^*(T^*M))\}.\] Furthermore the following proposition holds
%\[L(\partial \Sigma)=\{\mbox{Lie algebroid morphisms between  } T([0,1]) \mbox{  and  } T^*M .\}\ proven that $L(\partial \Sigma)$ is a coisotropic submanifold of $\Phi_{\partial}$ \footnote{see proof in the Appendix.} and in \cite{Ca} the following Theorem is proven:
\begin{proposition}\label{Coiso} \cite{Ca}. 
The space $C_{\Pi}$ is a coisotropic submanifold of $\Phi_{\partial}$.
\end{proposition}
In fact, the converse of this proposition also holds in the following sense. If we define $S(X,\eta)$ and $C_{\Pi}$ 
in the same way as before, without assuming that $\Pi$ satisfies Equation
(\ref{SN}) it can be proven that

\begin{proposition}\cite{Coi} \cite{Relational}.
If $C_{\Pi}$ is a coisotropic submanifold of $\Phi_{\partial}$, then $\Pi$ is a Poisson bivector field.
\end{proposition}
The geometric interpretation of the Poisson sigma model will lead us to the connection between Lie algebroids and Lie groupoids in Poisson geometry.
First we need some definitions.\\
A pair $(A, \rho)$, where $A$ is a vector bundle over $M$ and $\rho$ (called the anchor map) is a vector bundle morphism from $A$ to $TM$ is called a 
\textit{Lie algebroid} if
\begin{enumerate}
 \item There is Lie bracket $[,]_{A}$ on $\Gamma(A)$ such that the induced map $\rho_{*}: \Gamma (A) \to \mathfrak{X}(M)$ is a Lie algebra homomorphism.
\item \textit{Leibniz identity:} \[[X,fY]_{A}=f[X,Y] + \rho_{*}(X)(f)Y, \forall\, X,Y \in \Gamma (A), f \in \mathcal{C}^{\infty}(M).\]
\end{enumerate}
Lie algebras, Lie algebra bundles and tangent bundles appear as natural examples of Lie algebroids. For our purpose, the cotangent bundle of a 
Poisson manifold $T^*M$, where $[,]_{T^*M}$ is the Koszul bracket for 1-forms, that is defined 
for exact forms by
$$[df,dg]:=d\{f,g\}, \forall f,g \in \mathcal C ^{\infty}(M),$$

whereas for general forms it is recovered by Leibnitz and the anchor map given by
$\Pi^{\#}: T^*M \to TM$, is a central example of Lie algebroids.
To define a morphism of Lie algebroids we consider the complex $\Lambda^{\bullet}A^{*}$, where $A^*$ is the dual bundle and a differential $\delta_A$ is defined by the rules
\begin{enumerate}
 \item \[\delta_A f:= \rho^* df,\,\forall f \in \mathcal{C}^{\infty}(M).\] 
 \item \[\langle \delta_A \alpha, X \wedge Y \rangle:= -\langle \alpha,\, [X,Y]_A  \rangle + \langle \delta \langle \alpha, X \rangle, Y \rangle -\langle \delta \langle \alpha, Y \rangle, X \rangle, \, \forall X,Y \in \Gamma(A), \alpha \in \Gamma (A^*),\]
where $\langle, \rangle$ is the natural pairing between $\Gamma(A)$ and $\Gamma(A^*)$.
\end{enumerate}
A vector bundle morphism $\varphi: A \to B$ is a Lie algebroid morphism if 
\[\delta_A \varphi^*= \varphi^* \delta_B.\]
This condition written down in local coordinates gives rise to some PDE's the anchor maps and the structure functions for $\gamma(A)$ and 
$\Gamma(B)$ should satisfy. In particular, for the case of Poisson manifolds, $C_{\Pi}$ corresponds to the space of Lie algebroid morphisms between  $T (\partial \Sigma)$ and  $T^*M$
where the Lie algebroid structure on the left is given by the Lie bracket of vector fields on $T (\partial \Sigma)$ with 
identity anchor map and on the right is the one induced by the Poisson structure on $M$.

As it was mentioned before, it can be proven that this space, is a coisotropic submanifold of $T^*PM$. Its symplectic reduction, i.e. 
the space of leaves of its characteristic foliation, called the reduced phase space of the PSM, when is smooth, has a particular feature, 
it is a symplectic groupoid over $M$ \cite{Ca}. More precisely, a groupoid is a small category with invertible morphisms. When the spaces of objects and morphisms are 
smooth manifolds, a Lie groupoid over $M$, denoted by $G \rightrightarrows M$, can be rephrased as the following data\footnote{$G\times_{(s,t)}G$ is a smooth manifold whenever $s$ (or $t$) is a surjective submersion.}
\\
\xymatrixrowsep{4pc} \xymatrixcolsep{3pc} \xymatrix{
    &\,\,\,\,\;\;\,\,\;\,G\times_{(s,t)}   G  \ar[r]^{\,\,\,\,\,\;\;\;\;\;\;\;\;\mu}  & G \ar[r]^{i}   &G \ar@/_/[r]_t  \ar@/^/[r]^s & M \ar[l]_{\varepsilon}  & 
    }
\\
where $s,\,t,\, \iota,\, \mu$ and $\varepsilon$ denote the source, target, inverse, multiplication and unit map respectively,
such that the following axioms hold (denoting $G_{(x,y)}:= s^{-1}(x) \cap t^{-1}(y)$):
\\
    \textbf{(A.1)} $s\circ \varepsilon= t\circ \varepsilon =id_M$\\
    \textbf{(A.2)} If $g \in G_{(x,y)}$ and $h \in G_{(y,z)}$ then $\mu (g,h)\in G_{(x,z)}$\\
    \textbf{(A.3)} $\mu(\varepsilon\circ s \times id_{G})=\mu(id_{G}\times \varepsilon\circ t)=id_{G}  $ \\
    \textbf{(A.4)} $\mu(id_{G} \times i)=\varepsilon \circ t$\\
    \textbf{(A.5) }$\mu(i \times id_{G})=\varepsilon \circ s $\\
    \textbf{(A.6)} $\mu(\mu \times id_{G})=\mu (id_{G} \times \mu)$.\\
A Lie groupoid is called \textit{symplectic} if there exists a symplectic structure $\omega$ on $G$ such that
\[Gr_{\mu}:=\{ 
(a,b,c) \in G^3 \vert\,  c= \mu(a,b)\}\]
is a lagrangian submanifold of $G\times G \times \overline{G}$, where $\overline{G}$ denotes the sign reversed symplectic strucure on $G$. Finally, we can state the following

\begin{theorem} \label{Cafel}\cite{Ca}. The symplectic reduction $\underline{C_{\Pi}}$ of $C_\Pi$ (the space of leaves of the characteristic foliation), if it is smooth,
is a symplectic groupoid over $M$. 

%symplectic Lie groupoid integrating the Poisson manifold $M$.

\end{theorem}

The smoothness of the reduced phase space has particular interest. In \cite{Cra}, the necessary and suficient conditions for integrability of Lie algebroids, i.e. whether a Lie groupoid such that its \textit{infinitesimal version} 
corresponds to a given Lie algebroid exists, are stated. In \cite{PoissonCra}, these conditions have been further specialized to the Poisson case. It turns out that the reduced phase space of the PSM coincides with the space of 
equivalent classes of what are called $\mathcal A-$paths modulo $\mathcal A-$ homotopy \cite{Cra}, with $\mathcal A = T^*M.$ 

\section{The version before reduction.}
The main motivation for introducing the relational groupoid construction is the following. In general, 
the leaf space of a characteristic foliation is not a smooth finite dimensional manifold 
and in this particular situation, the smoothness of the space of reduced boundary fields is controlled 
by the integrability conditions stated in \cite{Cra}. 
In this paper, we define a groupoid object in the extended symplectic category, where the objects are symplectic manifolds, 
possibly infinite dimensional, and the morphisms are immersed Lagrangian submanifolds.  
It is important to remark here that this is extended category is not properly a category! 
(The composition of morphisms is not smooth in general). 
However, for our construction, the corresponding morphisms will be composable. 

We restrict ourselves to the case when $\mathcal{C}$ is the sometimes called \textit{Extended Symplectic Category}, 
denoted by $Sym^{Ext}$ and defined as follows:
\begin{definition}
$\mbox{Sym}^{Ext}$ is a category in which the objects are symplectic manifolds and the morphisms are immersed canonical relations. 
\footnote{This is not exactly a category because the composition of canonical relations is not in general a smooth manifold } Recall that 
$L: \mathcal M \nrightarrow \mathcal N$ is an immersed canonical relation between two symplectic manifolds $\mathcal M$ and $\mathcal N$ by definition if $L$ is an immersed Lagrangian submanifold of 
$\bar{\mathcal M}\times \mathcal N.$ \footnote{Observe here that usually one considers embedded Lagrangian submanifolds, but we consider immersed ones.}
$\mbox{Sym}^{Ext}$ carries an involution $\dagger:(\mbox{Sym}^{Ext})^{op} \to  \mbox{Sym}^{Ext} $ that is the identity in objects and in morphisms, for $f: A \nrightarrow B$, $f^{\dagger}:= \{(b,a) \in B\times A \vert (a,b) \in f\}$.
\end{definition}
This category extends the usual symplectic category in the sense that the symplectomorphisms can be thought in terms of canonical relations.

\begin{definition}
A \textbf{relational symplectic groupoid} is a triple $(\mathcal G,\, L,\, I)$ where 
\begin{itemize}
 \item $\mathcal G$ is a symplectic manifold (possibly infinite dimensional, 
in this case is a weak symplectic manifold. \footnote{it means that the induced map $T\mathcal G \to T^* \mathcal G$ is injective.})
 \item $L$ is an immersed Lagrangian submaifold of $\mathcal G ^3.$
 \item $I$ is an antisymplectomorphism of $G$
\end{itemize}
satisfying the following axioms
\end{definition}
\begin{itemize}
 \item \textbf{\underline{A.1}} $L$ is cyclically symmetric (i.e. if $(x,y,z) \in L$, then $(y,z,x) \in L$)

\item \textbf{\underline{A.2}} $I$ is an involution (i.e. $I^2=Id$).\\
\textbf{\underline{Notation}} $L$ is a canonical relation $\mathcal G \times \mathcal G \nrightarrow \bar{\mathcal G}$ and will 
be denoted by $L_{rel}.$ Since the graph of $I$ is a Lagrangian submanifold of $\mathcal G \times G$, $I$ is a canonical 
relation $\bar G \nrightarrow G$ and will be denoted by $I_{rel}$.\\
$L$ and $I$ can be regarded as well as canonical relations 
\[ \bar {\mathcal G} \times \bar{\mathcal G} \nrightarrow \mathcal G \mbox{ and } \mathcal G \nrightarrow \bar{\mathcal G}\]
respectively and will be denoted by $\overline{L_{rel}}$ and $\overline{I_{rel}}.$ The transposition 
\begin{eqnarray*} 
T: \mathcal G \times \mathcal G &\to& \mathcal G \times \mathcal G\\
(x,y) &\mapsto& (y,x)
\end{eqnarray*}
induces canonical relations
\[T_{rel}: \mathcal G \times \mathcal G \nrightarrow \mathcal G \times \mathcal G \mbox{ and } \overline{T_{rel}}: 
\bar {\mathcal G}\times \bar{\mathcal G} \nrightarrow \bar {\mathcal G}\times \bar{\mathcal G}.\]
The identity map $Id: \mathcal G \to \mathcal G$ as a relation will be denoted by $Id_{rel}: \mathcal G \nrightarrow \mathcal G$ and by $\overline{Id_{rel}}: \overline{\mathcal G} \nrightarrow \overline{\mathcal G}. $
\item \textbf{\underline{A.3}} 
$I_{rel} \circ L_{rel}= \overline{L}_{rel} \circ \overline{T}_{rel}\circ (\overline{I_{rel}} \circ \overline{I_{rel}}): \mathcal G \times \mathcal G \nrightarrow \mathcal G.$\\
\textbf{Remark 1:} Since $I$ and $T$ are diffeomorphisms, both sides of the equaliy correspond to immersed Lagrangian submanifolds.\\

Define \[L_3:= I_{rel} \circ L_{rel}: \mathcal G \times \mathcal G \nrightarrow \mathcal G.\] As a corollary of the previous axioms we get that
\begin{corollary}
$\overline{I_{rel}}\circ L_3 = \overline{L_3} \circ \overline{T_{rel}}\circ(\overline{I_{rel}}\times \overline{I_{rel}}).$ 
\end{corollary}
\item \textbf{\underline{A.4}} $L_3 \circ (L_3 \times Id)= L_3\circ (Id \times L_3): \mathcal G ^3 \nrightarrow G$ is an immersed Lagrangian submanifold. \\

The fact that the composition is Lagrangian follows from the fact that, since $I$ is an antisymplectomorphism, its graph is Lagrangian, 
therefore $L_3$ is Lagrangian, and so $(Id \times L_3)$ and $(L_3 \times Id).$
The graph of the map $I$, as a relation $* \nrightarrow \mathcal G \times \mathcal G$ will be denoted by $L_I$.\\
\item \textbf{\underline{A.5}} $L_3 \circ L_I$ is an immersed Lagrangian submanifold of $\mathcal G$.\\
\textbf{Remark 2:} It can be proven that Lagrangianity in these cases is automatical if we start with a finite dimensional symplectic manifold $\mathcal G$.

Let $L_1:= L_3 \circ L_I: * \nrightarrow \mathcal G $. From the definitions above we get the following
\begin{corollary}
\[\overline{I_{rel}}\circ L_1= \overline{L_1},\]
that is equivalent to 
\[I(L_1)= \overline{L_1}\]
where $L_1$ is regarded as an immersed Lagrangian submanifold of $\mathcal G$.
\end{corollary}
\item \textbf{\underline{A.6}}  \[L_3\circ(L_1 \times L_1)= L_1.\]

%\textit{Proof: }

\item \textbf{\underline{A.7}} $L_3\circ (L_1 \times Id)$ is an immersed Lagrangian submanifold of $\overline{\mathcal G} \times \mathcal G.$\\
We define
$$L_2:=L_3\circ (L_1 \times Id): \mathcal G \nrightarrow \mathcal G.$$
\begin{corollary}
$$L_2:=L_3\circ (Id \times L_1).$$ 
\end{corollary}
%\textit{Proof: }
\begin{corollary} $L_2$ leaves invariant $L_1, \, L_2$ and $L_3$, i.e.
\begin{eqnarray*}
L_2\circ L_1&=& L_1\\
L_2\circ L_2&=& L_2\\
L_2\circ L_3&=& L_3.
\end{eqnarray*}
%\textit{Proof: }
 
\end{corollary}
%Denoting by $L_2^{rel}$ the submanifold $L_3\circ (Id \times L_1)$ as a relation $* \nrightarrow \mathcal G \times \mathcal G$ 
\begin{corollary}\label{conv}
$$\overline{I_{rel}}\circ L_2= \overline{L_2}\circ \overline{I_{rel}}\mbox{ and } L_2^*=L_2.$$

\end{corollary}
\end{itemize}
                                                                                                                                                                                                                                                                                                                                                                                                                                                                  
The next set of axioms defines a particular type of relational symplectic groupoids, in which the relation 
$L_2$ plays the role of an equivalence relation and it allows to study the case of symplectic reductions.
\begin{definition}
A relational symplectic groupoid $(\mathcal G,\, L,\, I)$ is called \textbf{regular} if the following axioms are satisfied. Consider $\mathcal G$ 
as a coisotropic relation $* \nrightarrow \mathcal G$ denoted by $\mathcal G_{rel}$.
\end{definition}
\begin{itemize}
 \item \textbf{\underline{A.8}} $L_2 \circ \mathcal G_{rel}$ is an immersed coisotropic relation. \\
 \textbf{Remark 3:} Again in this case, the fact that this is a coisotropic relation follows automatically in the finite dimensional setting.

\begin{corollary} Setting $C:= L_2 \circ \mathcal G_{rel}$ the following corollary holds.\\
\begin{enumerate}
\item \[C^*= \mathcal G ^* \circ L_2\] 
 \item $L_2$ defines an equivalence relation on $C$.
 \item This equivalence relation is the same as the one given by the characteristic foliation on $C$. 
\end{enumerate}
\end{corollary}
\end{itemize}
%\textit{Proof: }
\begin{itemize}
\item \textbf{\underline{A.9}} The reduction $\underline{L_1}= L_1 / L_2$ is a finite dimensional smooth manifold. 
We will denote $\underline{L_1}$ by $M$.
\item \textbf{\underline{A.10}} $S:= \{(c,[l]) \in C \times M: \exists l \in [l], g \in \mathcal G \vert (l,c,g) \in L_3\}$ is an immersed submanifold 
of $\mathcal G \times M $.

\begin{corollary}
\[T:= \{(c,[l]) \in C \times M: \exists l \in [l], g \in \mathcal G \vert (c,l,g) \in L_3\}\] is an immersed submanifold of $\mathcal G \times M$.
\end{corollary}
The following conjectures (this is part of work in progress) give a connection between the symplectic structure on $\mathcal G$ and Poisson structures on $M$.

\begin{conjecture} Let $(\mathcal G, L,I)$ be a regular relational symplectic groupoid. Then, 
\begin{enumerate}
 \item There exists a unique Poisson structure on $M$ such that $S$ is coisotropic
in $\mathcal G \times M.$
\item This is also the unique Poisson structure on $M$ such that $T$ is coisotropic in $\mathcal G \times M.$
\end{enumerate}
\end{conjecture}
\begin{conjecture}
Assume $G:= C / L_2$ is smooth. Then $G$ is a symplectic groupoid on $M$ with structure maps $s:= S/L_2,\, t:= T/L_2, \mu:= L_{rel}/L_2, \iota= I, \varepsilon = L_1/L_2.$
\end{conjecture}

\end{itemize}

\end{section}

%\end{section}

%Observe that the space of morphisms over one object has a group structure and are called the \textit{isotropy groups} of the groupoid.

%$\bar {\mathcal G} \times \bar \mathcal G \nrightarrow \mathcal G$ and  %$ \mathcal G \nrightarrow \bar \mathcal G$ respectively, will be denoted by $\bar{L_{rel}}$ and $\bar{I_{rel}}$
%\textbf{(R.A.1)} $L_2$ is an equivalence relation on a submanifold $C$ of $\mathcal G$ .\\

%\textbf{(R.A.4)} $L_3\circ (id_{\mathcal{G}}\times \iota)\circ \triangle_2 \subseteq \L_2  \circ \varepsilon \circ t $. \\
%\textbf{(R.A.5)} $L_3\circ (\iota \times id_{\mathcal{G}})\subseteq \L_2  \circ \varepsilon \circ s$\\
%\textbf{(R.A.6)} $s( L_2 \circ \varepsilon )=id_M$ and $(L_2 \circ \varepsilon )(M)$ is a canonical relation.\\
%\textbf{(R.A.7)} $L_2\circ L_1=L_1.$

%In the sequel, $L_1$ will denote the composition $L_2 \circ \varepsilon$ and will be understood as a canonical relation $L_1: \emptyset \to \mathcal{G}$.
%\begin{center}
%\begin{figure}
%\psfrag{Rq}{$\mathbb{R}^q$}
%\centering%
%\center{\includegraphics[scale=0.75]{RA2}}
%\caption{The relational axioms.}
%\label{fig:FigureExample}
%\end{figure}
%\end{center}
\begin{definition} 
A \textbf{morphism between relational symplectic  groupoids} $(G,L_G,I_G)$ and $(H,L_H,I_H)$ 
is a map $F$ from $G$ and $H$ satisfying the following properties:
\begin{enumerate}
 \item $F$ is a Lagrangian subspace of $G\times \bar{H}$.
 \item $F\circ I_G=I_H \circ F$. 
 \item $F^3(L_G)=L_H.$
\end{enumerate}
\end{definition}
\begin{definition}

A morphism of relational symplectic groupoids $F: G \to H $ is called an \textbf{equivalence} 
if the transpose canonical relation $F^{op}$ is also a morphism.

\end{definition}
\textbf{Remark 4:} For our motivational example, it can be proven that
\begin{enumerate}
\item Different differentiability degrees (the $\mathcal C^{k}-$ type of the maps $X$ and $\eta$) give raise to equivalent relational symplectic groupoids.
\item For regular relational symplectic groupoids, $\mathcal G$ and $G$ are equivalent.
\end{enumerate}

\subsection{Examples}
The following are natural examples of relational symplectic groupoids.

%\begin{example}
 
%Let $(G,\omega)$ be a symplectic manifold with $\omega$ a multiplicative form. A relational symplectic groupoid over $G$ is defined as follows:
%\begin{eqnarray*}
%\mathcal{G}&=& G\times G.  \\
%M&=& G.\\
%L_1&=&\mathbb{R}^n \oplus \Omega^{1}_{0,ex}(\mathbb{R}^n).\\
%L_2&=& \triangle_2(G \times G) .\\
%L_3&=&\{(g_1,\,g_2,\,g_3,\,g_4,\,g_5,\,g_6) \in G^6\vert\, g_2=g_3,\, g_1=g_5,\, g_2=g_6\} . \\
%\varepsilon&=& \triangle_2: G \to G\times G.\\
%\iota&=& (g_1,g_2) \to (g_2,g_1).\\
%s&=& Pr_1.\\
%t&=&Pr_2.
%\end{eqnarray*}

%\end{example}

\subsubsection{Lie groupoids.}
Symplectic groupoids
Given a Lie symplectic groupoid $G$ over $M$, we can endow it naturally with a relational symplectic structure:
\begin{eqnarray*}
\mathcal{G}&=&G.\\
%L_1&=&\varepsilon(M).\\
%L_2&=&\triangle_2(G).\\ 
L&=&\{(g_1,g_2,g_3) \vert (g_1,g_2) \in G\times_{(s,t)}G,\, g_3=\mu(g_1,g_2)\}.\\
%S &=&\{(g,x)\in G \times M \vert s(g)=x \}.\\
%T &=&\{(g,x)\in G \times M \vert t(g)=x \}.\\
I &=& g \mapsto g^{-1},\, g \in G.
\end{eqnarray*}

\subsubsection{Symplectic manifolds with a given Lagrangian submanifold:} Let $(G,\omega)$ be a symplectic manifold and $\mathcal L$ a Lagrangian submanifold of $G$. We define
\begin{eqnarray*}
\mathcal{G}&=&G.\\
%M&=& pt.\\
%L_1&=&L.\\
%L_2&=&L\times L.\\ 
L&=&\mathcal L \times \mathcal L\times \mathcal L.\\
%S &=&\{(g,pt)\vert g\in G \}.\\
%T &=&\{(g,pt)\vert g\in G \}.\\
I&=&\{\mbox{identity of $G$}\}. 
\end{eqnarray*}
It is an easy check that this construction satisfies the relational axioms and furthermore
\begin{proposition} 
The previous relational symplectic groupoid is equivalent to the zero dimensional symplectic groupoid (a point with zero symplectic structure and empty relations).
\end{proposition}
\textit{Proof:} It is easy by checking that $L$ is an equivalence from the zero manifold to $\mathcal{G}$ .

 \subsubsection{Powers of symplectic groupoids:}

Let us denote $G_{(1)}=G$, $G_{(2)}$ the fiber product $G\times_{(s,t)}G$, $G_{(3)}=G \times_{(s,t)}(G \times_{(s,t)}G)$ and so on . It can be proven the following
\begin{lemma} \cite{Relational}. Let $G\rightrightarrows M$ be a symplectic groupoid.
\begin{enumerate}
 \item $G_{(n)}$ is a coisotropic submanifold of $G^n$.
\item The reduced spaces $\underline{G_{(n)}}$ are symplectomorphic to $G$. Furthermore, there exists a natural symplectic groupoid structure on $\underline{G_{(n)}}$ coming from the quotient, isomorphic to the groupoid structure on $G$. 
\end{enumerate}

\end{lemma}

We have natural canonical relations 
$P_n:\underline{G_{(n)}}\to G^n$ defined as:
\[P:=\{(x,\alpha, \beta)\vert x \in \underline{G_{(n)}},\, [\alpha]=[\beta]=x\},\]
 satisfying the following relations:
\[P^{op}\circ P= Gr(Id_G),\, P\circ P^{op}= \{(g,h)\in G^n \vert [g]=[h]\}.\] 
It can be checked that 
\begin{proposition}
$G_{(i)}$ is equivalent to $G_{(j)}, \forall i,\, j \geq 1$ and the equivalence
 is given by $P_i \circ P_j^{op}$.
\end{proposition}

\subsubsection{The cotangent bundle of the path space of a Poisson manifold.}

This is the motivational example for the construction of relational symplectic groupoids. In this case, the coisotropic submanifold 
$C_{\Pi}$ is equipped with an equivalence relation, called $T^*M$- homotopy \cite{Cra}, and denoted by $\sim$. More precisely, to points of $C_{\Pi}$ are $\sim$- equivalent if they belong to the same leaf of the characteristic foliation of $C_{\Pi}$.
We get the following relational symplectic groupoid (where $L$ is the restriction to the boundary of the solutions of the Euler-Lagrange equations in the bulk)
\begin{eqnarray*}
\mathcal{G}&=& T^*(PM).\\
%M&=& pt.\\
%L_1&=&L.\\
%L_2&=&L\times L.\\ 
L&=& \{(X_1,\eta_1),(X_2,\eta_2),(X_3,\eta_3) \in C_{\Pi}^3 \vert (X_1 *X_2, \eta_1* \eta_2)\sim (X_3 * \eta_3)\}.\\
%S &=&\{(g,pt)\vert g\in G \}.\\
%T &=&\{(g,pt)\vert g\in G \}.\\
I&=& (X,\eta)\mapsto (\phi^*X, \phi^*\eta)\}.
\end{eqnarray*}

Here $*$ denotes path concatenation and
\begin{eqnarray*}
\phi: [0,1]&\to&[0,1]\\
t&\mapsto& 1-t 
\end{eqnarray*}

\begin{theorem}\cite{Relational}. The relational symplectic groupoid $\mathcal G$ defined above is regular.

\end{theorem}
The improvement of Theorem \ref{Cafel} in terms of the relational  symplectic groupoids can be summarized as follows.  $L_1$ can be understood as the space of $T^*M$- 
paths that are $T^*M$- homotopy equivalent to the trivial $T^*M$- paths and 
$$\overline{L_1}:= \cup_{x_0 \in M} T^*_{(\overline{X,\eta})}PM \cap L_1,$$
where $(\overline{X,\eta})= \{(X,\eta) \vert X \equiv X_0, \eta \in \ker \Pi^{\#}\},$ we can prove the following

\begin{theorem}\cite{Relational}. \textit{If the Poisson manifold $M$ is integrable, then, there exists a tubular neigborhood of the zero section of $T^*PM$, denoted by $N(\Gamma_0(T^*PM))$ such that $\overline{L_1}\cap N(\Gamma_0(T^*PM))$ is an embedded submanifold of $T^*PM$.} \\
\end{theorem}
\begin{theorem}\cite{Relational}. \textit{If M is integrable, then $L_1\cap N(\Gamma_0(T^*PM)),\, L_2\cap N(\Gamma_0(T^*PM))^2$ and $L_3\cap N(\Gamma_0(T^*PM))^3$ are embedded Lagrangian submanifolds.} \\
\end{theorem}


\begin{thebibliography}{10}

\bibitem{Duistermaat} J.J. Duistermaat and J.A.C. Kolk, \emph{Lie Groups} Universitext, Springer, 1999.
\bibitem{Ca} A. S. Cattaneo and G. Felder. \emph{Poisson sigma models and symplectic groupoids}, Progress in Mathematics 198, 61-93, 2001.



%\bibitem{Frobenius} J. Teichmann. \emph{A Frobenius Theorem on Covenient Manifolds.}

\bibitem{Coiso} A.S.  Cattaneo. \emph{Coisotropic submanifolds and dual pairs.}, preprint.

%\bibitem{Michor1} P. Michor. \emph{Manifolds of Differentiable Mappings,} Shiva Mathematics Series, 1980.

%\bibitem{Milnor} J.Milnor. \emph{Morse Theory}, Annals of Mathematics Studies, No. 51, Princeton University Press, 1963.

%\bibitem{Igusa} K. Igusa. \emph{Iterated integrals and superconnections}, arXiv:0912.0249v1, 2009.

%\bibitem{PoissonStr} J.P. Dufour and N.T. Zung. \emph{Poisson Structures and Their Normal Forms}, Progress in Mathematics, Birkheauser Verlag, 2005.

%\bibitem{Montreal} A. Cattaneo, \emph{Integration of Poisson manifolds} April 2010, Montreal, Private Communication.

%\bibitem{CattaneoDeformation} A. Cattaneo, \emph{Deformation quantization and reduction}, Contemporary Mathematics, 450: 79-101. 2008.

%\bibitem{Brown} R. Brown and O. Mucuk, \emph{Foliations, locally Lie groupoids and holonomy}

%\bibitem{Marle} C.M. Marle, \emph{From momentum maps and dual pairs to symplectic and Poisson groupoids} ? 

\bibitem{Branes} A. S. Cattaneo and G.Felder \emph{Coisotropic submanifolds in Poisson geometry and branes in the Poisson sigma model}, Lett. Math. Phys. 69. 157-175. 2004.

%\bibitem{Dazord} A. Coste, P. Dazord, and A Weinstein, \emph{Groupo\"ides Symplectiques} Publ. Dept. Math. Univ. Claude-Bernard Lyon I, 1987.

%\bibitem{Weinstein} A. Weinstein, \emph{Symplectic groupoids and Poisson manifolds} Bull. Amer. Math. Soc. 16 (1987), 101-104.

%\bibitem{Cannas} A. Cannas da Silva and A. Weinstein, \emph{Geometric models for noncommutative algebras}, American Mathematical Society, Providence, RI, 1999.

%\bibitem{Bursztyn} H. Bursztyn, \emph{A brief introduction to Dirac manifolds}, To appear in Geometric and Topological Methods for Quantum Field Theory, CUP, 2012.

\bibitem{Relational} A.S. Cattaneo and I. Contreras, \emph{Relational symplectic groupoids and Poisson sigma models with boundary}, in preparation.

%\bibitem{Heunen} C. Heunen, I. Contreras and A. Cattaneo, \emph{Relative Frobenius algebras are groupoids} http://arxiv.org/abs/1112.1284, to appear in J. of Pure and App. Algebra, 2012.

\bibitem{Coi} A. S. Cattaneo and G. Felder, \emph{ Coisotropic submanifolds and dual pairs}, unpublished.

\bibitem{Algebroid} A. Cattaneo, \emph{On the integration of Poisson manifolds, Lie algebroids, and Coisotropic submanifolds}, Letters in Mathematical Physics 67: 33-48, 2004.

\bibitem{AlbertoPavel1}A. S. Cattaneo, P. Mnev and N. Reshetikhin, \emph{Classical BV theories on manifolds with boundaries}. math-ph/1201.0290
\bibitem{AlbertoPavel2}A. S. Cattaneo, P. Mnev and N. Reshetikhin, \emph{Classical and Quantum Lagrangian Field Theories with Boundary}, preprint. 
\bibitem{Chenchang} H.H Tseng, C. Zhu, \emph{Integrating Lie algebroids via stacks},
Compositio Mathematica, Volume 142 (2006), Issue 01, pp 251-270.

\bibitem{Corfu} A.S. Cattaneo, P. Mnev and N. Reshetikin, \emph{Classical and Quantum  Lagrangian Field Theories with Boundary}, Proceedings of the Corfu Summer Institute 2011 School and Workshops on Elementary Particle Physics
and Gravity, Corfu, Greece, 2011.
\bibitem{Cra} M. Crainic and R. L. Fernandes. \emph{ Integrability of Lie brackets}. Ann. of Math.(2) 157 (2003), 575--620.
%\bibitem{Cra} M. Crainic and R. L. Fernandes.  \emph{ Integrability of Lie brackets}.(2) 157 (2003), 575?620.


\bibitem {PoissonCra}M. Crainic and R. L. Fernandes. \emph{ Integrability of Poisson brackets},  parXiv:math/0210152v1

\bibitem{Mor} I. Moerdjik and J.Mrcun. \emph{Intoduction to Foliations and Lie Groupoids}, Cambridge studies in advanced mathematics, 91, 2003.
\bibitem{Severa} P. Severa, \emph {Some title containing the words homotopy and symplectic, e.g. this one}, preprint math.SG-0105080.



\end{thebibliography}
\end{document}